\documentclass{SIG-ALTERNATE} 

\newtheorem{thm}{Theorem}[section]
\newtheorem{prop}[thm]{Proposition}

\usepackage{bbm}
\newcommand{\ind}[1]{\mathbbm{1}_{\{#1\}}}   
\newcommand{\E}[1]{\mathbb{E}[#1]}

\newcommand{\Prob}{\mathbb{P}}
\newcommand{\htt}{\hat}

\newcommand{\blt}{\boldsymbol}

\newcommand{\wed}{\wedge}

\newcommand{\be}{\beta}
\newcommand{\R}{\mathbb{R}}
\graphicspath{{Figures/}}

\confinfo{MAMA 2017 Workshop Urbana-Champaign, Illinois, USA.}
\title{Electric vehicle charging: a queueing   approach}

\numberofauthors{3}
\author{
\alignauthor
Angelos Aveklouris\\
       \affaddr{Eindhoven University of Technology}\\
       \email{a.aveklouris@tue.nl}
\alignauthor
Yorie Nakahira\\
       \affaddr{California Institute of Technology }\\
       \email{ynakahir@caltech.edu}
\alignauthor Maria Vlasiou\\
       \affaddr{Eindhoven University of Technology}\\
       \email{m.vlasiou@tue.nl}
\and  
\alignauthor Bert Zwart\\
       \affaddr{Eindhoven University of Technology,}\\
       \affaddr{Centrum Wiskunde en Informatica}\\
       \email{bert.zwart@cwi.nl}}

\begin{document}
\maketitle
\begin{abstract}
The number of electric vehicles (EVs) is expected to increase. As a consequence, more EVs will need charging, potentially causing not only congestion at charging stations, but also in the distribution grid.
Our goal is to illustrate how this gives rise to resource allocation and performance problems that are of interest to the Sigmetrics community.
\end{abstract}

\section{Introduction}
EVs consume a large amount of energy and as a result the charging of EVs is causing congestion in the distribution grid \cite{clement2010impact}, which is exacerbated as the number of charging stations is limited. Motivated by this, we consider a stylized model that models the interplay of two sources of congestion (as not all cars find a space): (i) the number of available spaces with charging stations; (ii) the amount of available power.

Despite being a relatively new topic, the engineering literature on EV charging is huge. Here, we give only a sample. In \cite{su2012performance}, an algorithm for optimally managing a large number of plug-in EVs charging at a parking station is suggested. In \cite{sortomme2011coordinated}, optimal charging algorithms that minimize the impact of plug-in EV charging on the connected distribution grid are proposed. Examples of studies where randomness is taken into account are \cite{li2012modeling}, in which a methodology of modeling the overall charging demand of plug-in EVs is proposed, and \cite{turitsyn2010robust} where control algorithms based on randomized EV charging arrival time are suggested. Mathematical models where vehicles communicate beforehand with the grid to convey information about their charging status are studied in \cite{said2013queuing}. In \cite{kempker2016optimization}, cars are the central object and a dynamic program is formulated that prescribes how cars should charge their battery using price signals. Though the class of problem at hand fits well to the performance analysis, the only other line of work where such ideas are used is \cite{bayram2013electric} and \cite{yudovina2015socially}, where a gradient scheduler is proposed to minimize delays.

A common feature of the above studies is that they apply to shorter operational time-scales. Since the desired scale of infrastructure does not exist yet, it is important to consider models that can be used on longer (investment) time-scales. Equilibrium models are quite popular for investment and policy analysis of energy systems \cite{ehrenmann2011generation}.
We therefore consider a stylized equilibrium queueing model that takes into account both congestion in the distribution grid, as well as congestion in the number of available spaces with charging stations.
We consider a stylized model of a parking lot with finitely many spaces in which EVs (customers) arrive randomly in order to get charged (for another application of queueing theory to parking lots see \cite{larson2007congestion}). The EVs have a random parking time and a random energy demand.
Thus, each EV receives two kinds of service, parking and charging. We assume that all available power is charged at the same rate to all cars that need charging; some of our results can be extended
to time-varying arrival rates and multiple types of users and stations, but due to space we do not do so here.

Under Markovian assumptions, our analysis focuses on the probability that an EV leaves the parking lot with a fully charged battery. Specifically, we develop bounds and a fluid approximation, and report partial results on a diffusion approximation. Our mathematical results are closely related to work on processor-sharing queues with impatience \cite{gromoll2006impact}, though the model here is more complicated as there is limited number of spaces in the system and fully charged cars may not leave immediately as they are still parked.

\section{Model description}\label{S:model desctiption}
We consider a charging station with $K>0$ parking spaces, each having an EV charger.
We assume that the arrival, parking and charging times of EVs are mutually independent, and exponential with rates $\lambda, \mu$ and $\nu$, respectively.
EVs leaves the system after their parking time expires. An EV may leave the system without its battery being fully charged. Furthermore, if all spaces are occupied, a newly arriving EV does not enter the system but leaves immediately. As such, the total number of vehicles in the system can be modeled by an Erlang loss system, though we need a more detailed description of the state space.

We denote by $Q(t)\in \{0,1,\ldots, K\}$ the total number of EVs in the system at time $t\geq 0$, where $Q(0)$ is the initial number of EVs. Further, we denote by $U(t)\in\{0,1,\ldots, Q(t)\}$ the number of EVs of which their battery is not fully charged at time $t$ and by $U(0)$ the number of vehicles initially in the system. Thus, $C(t)=Q(t)-U(t)$ represents the number of EVs with a fully charged battery at time $t$.

The power consumed by the parking lot is limited and depends on the number of uncharged EVs at time $t$. We let it be given by
$f:\R_{+} \rightarrow \R_{+} $,
$f(U(t)):= \min\{U(t),M\}$.
 Here, $0<M\leq  K$ denotes the maximum number of cars the parking lot can charge at full power.

\section{MAIN RESULTS} \label{S:main results}
We present bounds and  approximations based on fluid and diffusion limits for the fraction of EVs that get fully charged. Proofs (and results for other performance measures) will be presented in an extended version of this paper.

\subsection{Bounds}
Under our assumptions, the number of uncharged EVs and the total number of EVs in the system $(U(t),Q(t))$, is a two dimensional Markov process.
The fraction of fully charged EVs in stationarity is given by the ratio: $\frac{\E{ C(\infty)}}{\E {Q(\infty)}}$. In the special case $K=M$, we can compute explicitly the joint distribution, and in the case  $K=\infty$, the distribution of the number of uncharged EVs is given by a variation of the Erlang A formula (see \cite{zeltyn2005call} for details on the Erlang A model). Note that, in our model customers / EVs can leave the system also during their service, unlike in the Erlang A queue.
Based on these two special cases ($K=M$, $K=\infty$), the following proposition, which can be proved using Markov-rewards methods, presents an upper and a lower bound for the fraction of EVs that get fully charged.
\begin{prop}
Let $C^{K}_M(\infty)$ and $Q^{K}_M(\infty)$ be the number of fully charged EVs and the total number of EVs in stationarity for the system $(K,M)$. We have that
\begin{equation}\label{In:bounds}
\frac{\E {C^{\infty}_M(\infty)}}{\E {Q^{\infty}_M(\infty)}} \leq
\frac{\E{ C^{K}_M(\infty)}}{\E {Q^{K}_M(\infty)}}\leq  \frac{\E {C^{K}_K(\infty)}}
{\E {Q^{K}_M(\infty)}}.
\end{equation}
\end{prop}

\subsection{Fluid approximation}
We develop a fluid approximation for finite $K$, following a similar approach as in \cite{gromoll2006impact}. The main differences are the finitely many servers in the system and that
the state space consists of two regions: $U(\infty)>M$ and $U(\infty)\leq  M$.

Consider a family of models as defined earlier indexed by $n$. The fluid scaling (in steady state) is given by $\frac{U_n(\infty)}{n}$. To obtain a non-trivial  fluid limit, we assume that the capacity of power in the $n^{\text{th}}$ system is given by $nM$, the arrival rate  by $n\lambda$, and the number of parking spaces by $nK$.

\begin{prop}
Let $E_\mu$ and $E_\nu$ be exponential random variables with rates $\mu$ and $\nu$. We have that $\frac{U_n(\infty)}{n} \rightarrow u^*$, as $n \rightarrow \infty$. In addition, $u^*$ is given by the unique positive solution of the following fixed-point equation:
\begin{equation*}\label{eq:fluid proxy}
u^* = \min\{\lambda, \mu K\} \E{ \min \{E_\mu , E_\nu \max\{ 1, \frac{ u^*}{M}\}\}}.
\end{equation*}
\end{prop}
Observe that if we define $f(U(\cdot))=1$ (i.e., the processor sharing discipline) and replace $K$ by $n \lambda K$ (assuming for simplicity $\mu=1$), we derive \cite[Equation 4.1]{gromoll2006impact}.

We directly use a modified form of our fluid approximation, which can be derived heuristically using Little's law and a version of the snapshot principle (essentially assuming an EV sees the system in stationarity throughout its sojourn). Let $P_K$ be the blocking probability in a loss system with $K$ servers.
 To obtain our approximation, we replace $\min \{\lambda, \mu K\}$ by $\lambda (1-P_K)$, leading to
\begin{equation}\label{eq:fluid proxy}
u^* = \lambda (1-P_K) \E{ \min \{E_\mu , E_\nu \max\{ 1, \frac{ u^*}{M}\}\}}.
\end{equation}

Let $P_s$ denote the probability that an EV leaves the parking lot with fully charged battery in the fluid model. It is given by
$P_s=\Prob(E_\mu > E_\nu \max\{ 1, u^*/M\})$, where $u^*$ is the unique solution of \eqref{eq:fluid proxy}. Under our assumptions, the explicit expression for this probability can be found. That is,
\begin{equation}\label{eq:fluidFrac}
P_s=
\begin{cases}
\frac{\nu}{\mu+\nu},
& \mbox{$u^*\leq M$}, \\
\frac{\nu M}{\lambda (1-P_K)},
& \mbox{$u^*> M$}.
\end{cases}
\end{equation}

\subsection{Diffusion Approximation}

Let $\beta$ and $ \kappa$ be real numbers. Consider the following asymptotic regime. Define
$M_n=\frac{\lambda_n}{\nu+\mu}+\beta \sqrt{n}$ and $\lambda_n=n(\nu+\mu)$, i.e., the ``square-root staffing rule'' as in \cite{fleming1994heavy} and
\cite{garnett2002designing}. In addition, define
$K_n=\frac{\lambda_n}{\mu}+\kappa \sqrt{n}$.
The diffusion scaling is given by
$\htt U_n(t):=\frac{U_n(t)-\frac{\lambda_n}{\nu+\mu}}
{\sqrt{n}}$
and
$\htt Q_n(t):=\frac{Q_n(t)-\frac{\lambda_n}{\mu}}{\sqrt{n}}$.
\begin{thm}
If $(\htt U_n(0),\htt Q_n(0)) \overset{d} \rightarrow
(\htt U(0),\htt Q(0))$ then  \\
$(\htt U_n(\cdot),\htt Q_n(\cdot)) \overset{d} \rightarrow
(\htt U(\cdot),\htt Q(\cdot))$,
as $n \rightarrow \infty$.
 The diffusion limit satisfies the following 2-dimensional stochastic differential equation
\begin{equation}\label{eq:SDE}
\begin{split}
\begin{bmatrix}
d\htt U(t) \\
d\htt Q(t)
\end{bmatrix}
=&
\begin{bmatrix}
\sqrt{2 (\nu+\mu)}& 0  \\
0 & \sqrt{2 (\nu+\mu)}
\end{bmatrix}
\begin{bmatrix}
dW_{\htt U}(t) \\
dW_{\htt Q}(t)
\end{bmatrix}\\
&+
\begin{bmatrix}
b_1(\htt U(t),\htt Q(t)) \\
b_2(\htt U(t),\htt Q(t))
\end{bmatrix}
dt
-
\begin{bmatrix}
dY(t) \\
dY(t)
\end{bmatrix},
\end{split}
\end{equation}
where $b_1(x,y)=-\nu (x \wed \beta)-\mu x$ and $b_2(x,y)=-\mu y$. Further, $W_{\htt U}(t)$ and $W_{\htt Q}(t)$ are driftless, univariate Brownian motions such that
$2(\nu+\mu)\E{W_{\htt U}(t) W_{\htt Q}(t)}=
(\nu+2\mu) t$. In addition, $Y(\cdot)$ is the unique nondecreasing nonnegative process such that \eqref{eq:SDE} holds and
$\int_{0}^{\infty} \ind{{\htt Q}(t)<\kappa}dY(t)=0$.
\end{thm}

Note that $\htt Q(t)$ satisfies the known Erlang B diffusion \cite{pang2007martingale}. When $\kappa=\infty$ the system \eqref{eq:SDE} has an explicit invariant distribution. Take the vectors $\blt{m}_{-}=(0,0)$, $\blt{m}_{+}=(-\frac{\nu \be}{\mu},0)$ and the positive definite matrices
$\blt{\Sigma}_{-}=
\begin{bmatrix}
1 & \frac{2}{\nu+2\mu} \\
\frac{2}{\nu+2\mu} & \frac{\nu+\mu}{\mu}
\end{bmatrix}$
and
$\blt{\Sigma}_{+}=
\begin{bmatrix}
\frac{\nu+\mu}{\mu} & \frac{1}{\mu} \\
\frac{1}{\mu} & \frac{\nu+\mu}{\mu}
\end{bmatrix}$. Let
$f_{-}$ and $f_{+}$
be 2-dimensional normal pdfs with mean vectors $\blt{m}_{-}$, $\blt{m}_{+}$ and covariance matrices $\blt{\Sigma}_{-}$ and $\blt{\Sigma}_{+}$, respectively. In case $K=\infty$, we can show that the joint steady state pdf of the random vector
$(\htt U(\infty),\htt Q(\infty))$
can be written as
$$ \phi(x,y)=c_1f_{-}(x,y)
\ind{x \leq  \beta}+c_2f_{+}(x,y)
\ind{x > \beta},$$
where $c_1,c_2$ are given in \cite[Equations 3.9--3.10]{fleming1994heavy}.

\subsection{Numerical evaluation and discussion}
In Fig.\ 1--3 we depict the bounds in \eqref{In:bounds} and the fluid approximation in \eqref{eq:fluidFrac} for 3 cases: moderately, critically, and over-loaded. Further, we fix $\nu=\mu=1$. The vertical axes give the probability that an EV leaves the parking lot with fully charged battery (success probability) and the horizontal axes give the ratio $M / K$. The lower and the upper bounds seem to be tight for $M>0.7K$. Also, the lower bound is tight under light load, in the other cases the fluid approximation works well for $K=50$.

\begin{figure}[!h]
	\centering
		\includegraphics[width=.45\linewidth]
		{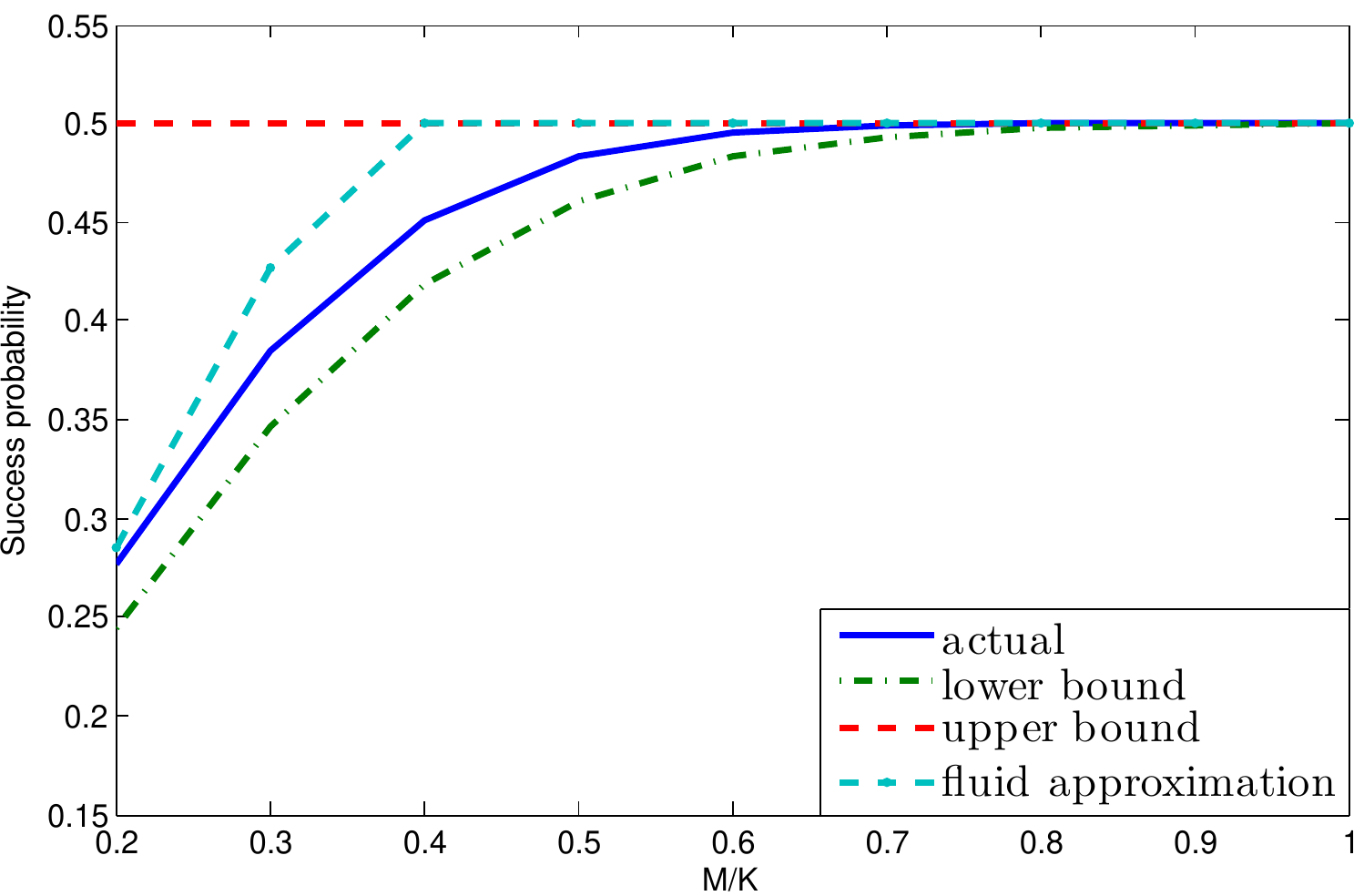}
		\includegraphics[width=.45\linewidth]
		{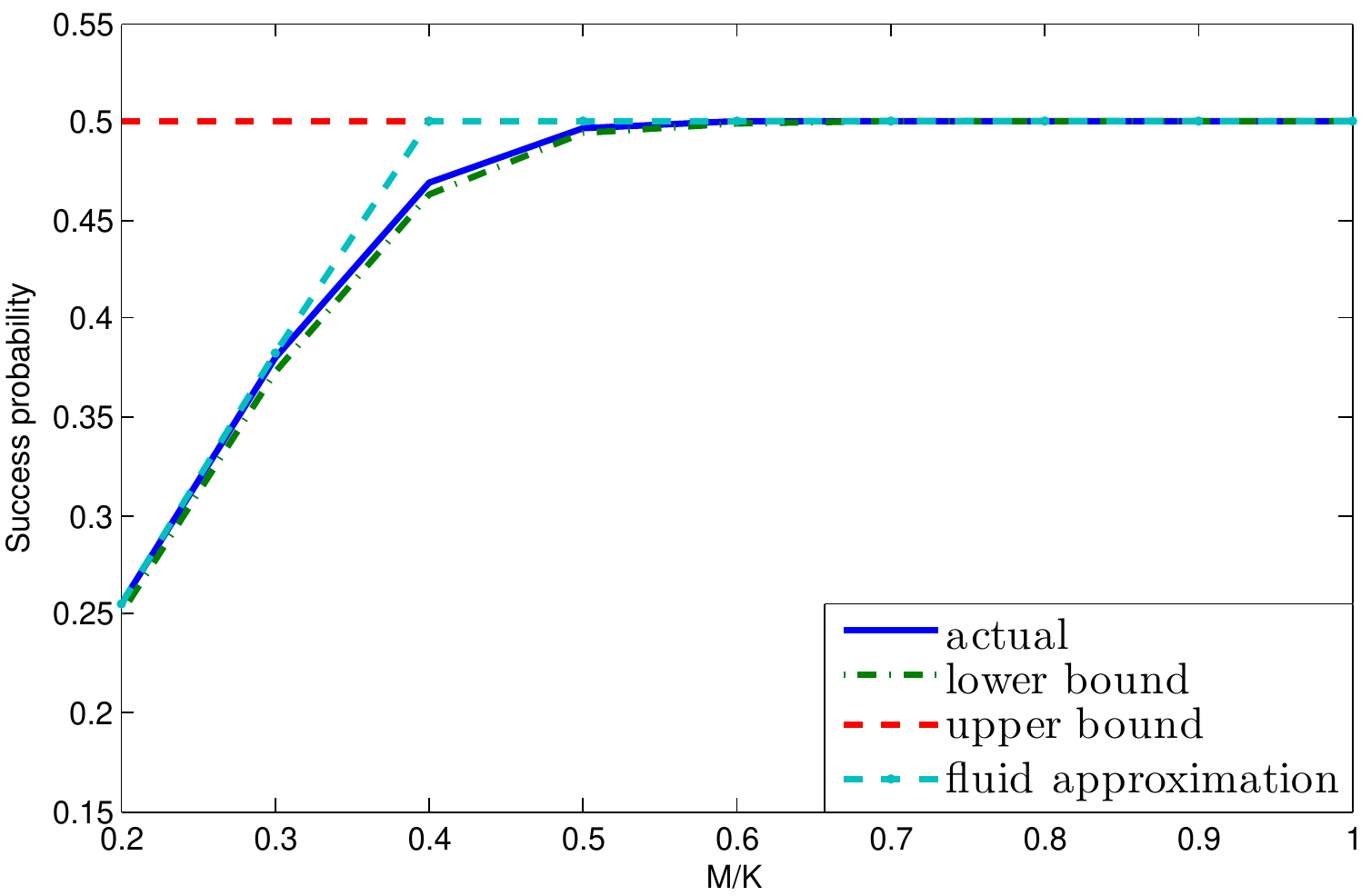}
	\caption{$K=10,50$ and $\lambda=0.8K$}
\end{figure}

\begin{figure}[h]
	\centering
	\includegraphics[width=.45\linewidth]
	{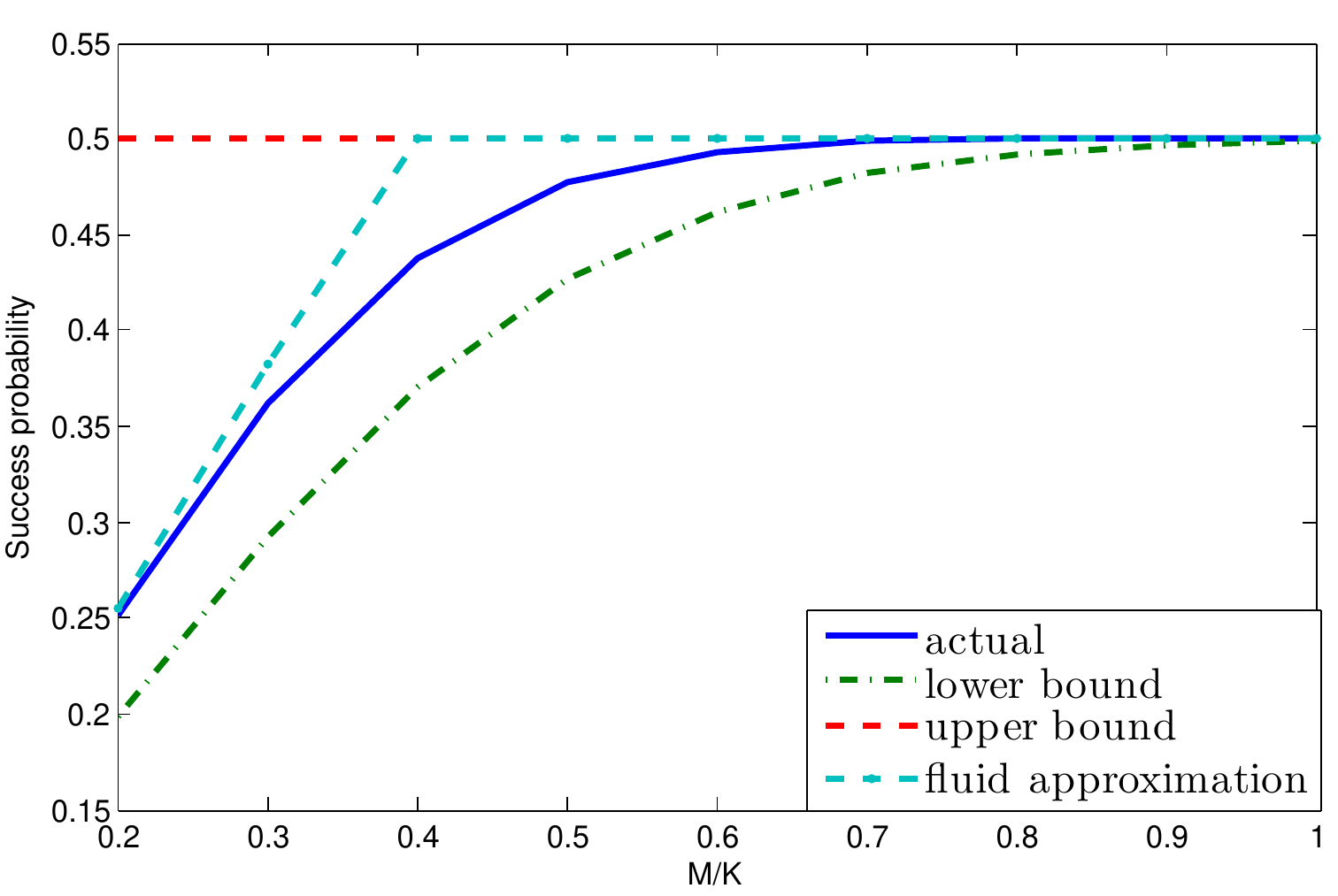}
	\includegraphics[width=.45\linewidth]
	{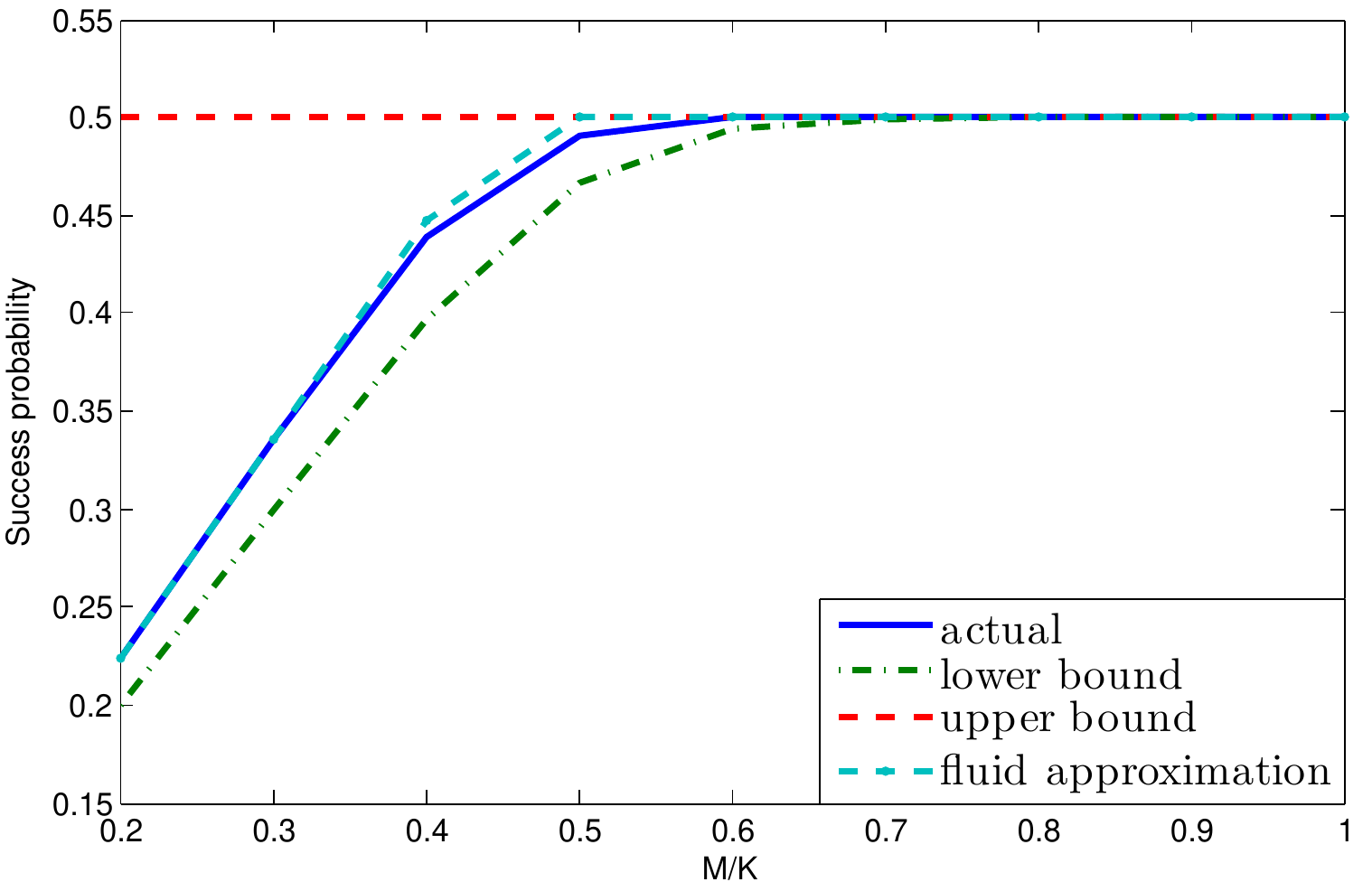}
	
	\caption{$K=10,50$ and $\lambda=K$}
\end{figure}

\begin{figure}[h!]
	\centering
		\includegraphics[width=.45\linewidth]
		{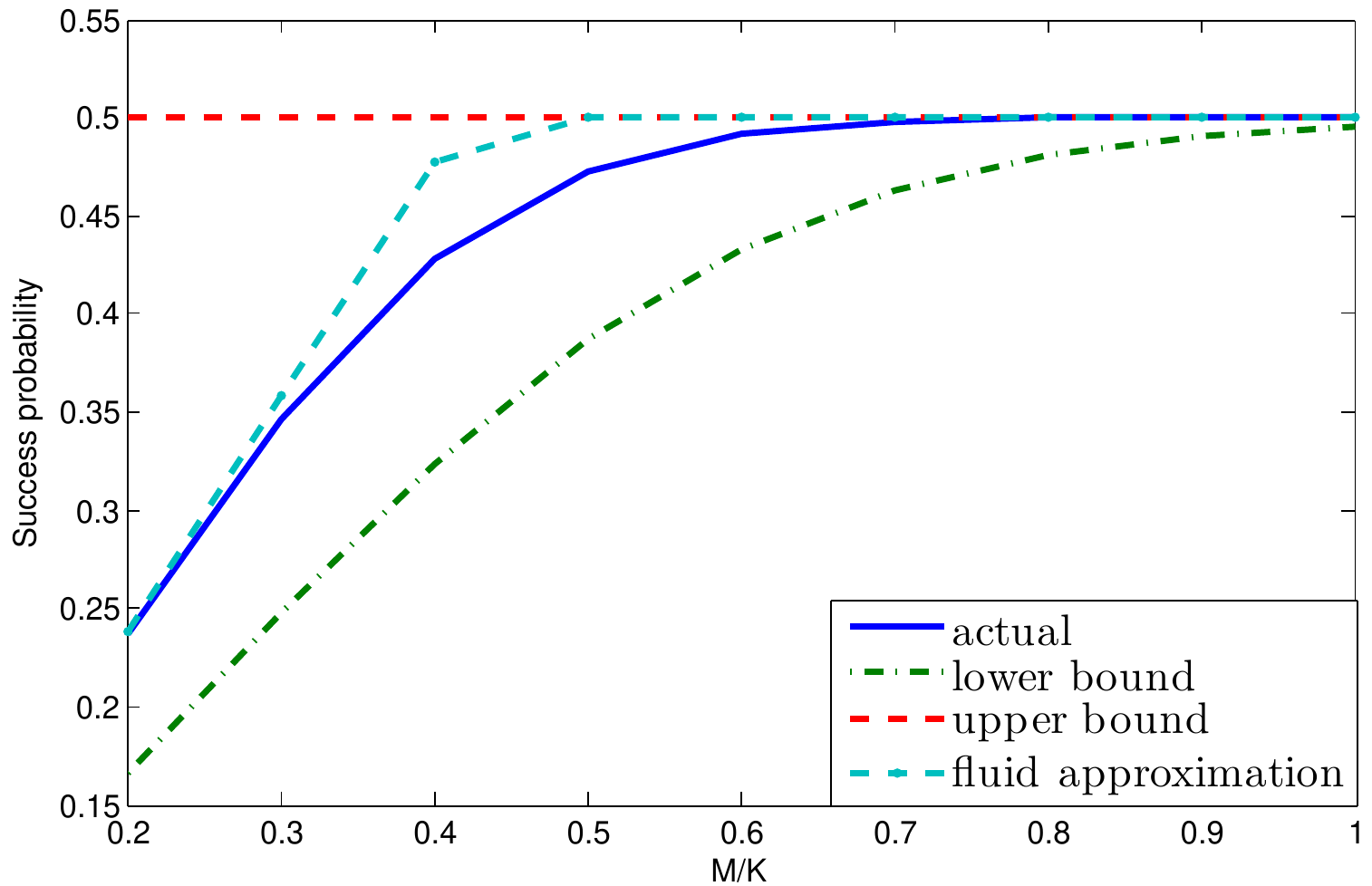}
		\includegraphics[width=.45\linewidth]
		{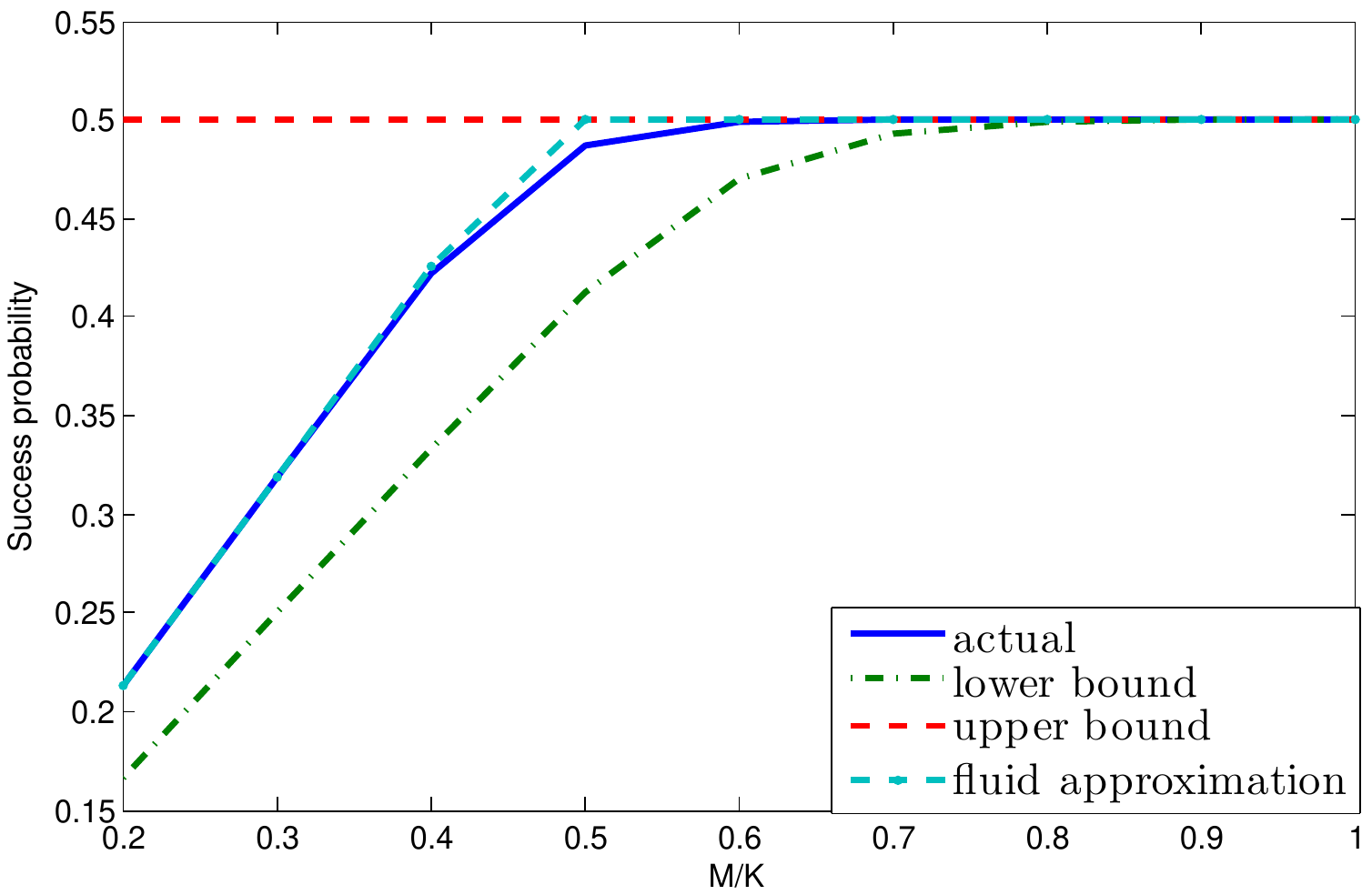}

	\caption{$K=10,50$ and $\lambda=1.2K$}
\end{figure}


\section{Discussion and extensions}
Our numerical results show there is room for improvement for critically loaded systems, making it worthwile to derive the invariant distribution of the process in  Theorem 3.4; the solution for $\kappa=\infty$ did not yield better results than the Erlang A lower bound.

From an applications standpoint, it is important to remove various model assumptions.
If parking and charging times are given by the (possibly dependent) generally distributed random variables $B$ and $D$, we can develop a measure-valued fluid model by extending  \cite{gromoll2006impact}. The fluid limit in steady-state will be defined by the  fixed point equation
\begin{equation*}
u^* = \lambda (1-P_K) \E{ \min \{B , D \max\{ 1, \frac{ u^*}{M}\}\}}.
\end{equation*}
We are currently extending this to time-varying arrival rates,
multiple customer classes, and multiple
parking lots. The distribution network for the latter is modeled as a tree network in \cite{carvalho2015critical} where simulation results are presented. On a high level, the analysis is reminiscent of \cite{remerova2014fluid}.

\bibliographystyle{abbrv}

\bibliography{sigproc}
\balancecolumns

\end{document}